\documentclass{amsart}
 \usepackage{epsfig}
\usepackage{graphicx}
\usepackage{amscd}
\usepackage{amsmath}
\usepackage{amsxtra}
\usepackage{amsfonts}
\usepackage{amssymb}

\newtheorem{theorem}{Theorem}[section]

\newtheorem{conjecture}[theorem]{Conjecture}
\theoremstyle{definition}
\newtheorem{definition}[theorem]{Definition}

\newtheorem{problem}{Problem}

\theoremstyle{remark}

\renewcommand{\theclaim}{\textup{\theclaim}}

\newtheorem*{acknowledgements}{Acknowledgements}

\numberwithin{equation}{section}

\def\openone%{\hbox{\upshape \small1\kern-3.3pt\normalsize1}}

{\mathchoice

{\hbox{\upshape \small1\kern-3.3pt\normalsize1}}

{\hbox{\upshape \small1\kern-3.3pt\normalsize1}}

{\hbox{\upshape \tiny1\kern-2.3pt\SMALL1}}

{\hbox{\upshape \Tiny1\kern-2pt\tiny1}}}

\makeatletter

\newbox\ipbox

\newcommand{\ip}[2]{\left\langle #1\, , \,#2\right\rangle}
\newcommand{\diracb}[1]{\left\langle #1\mathrel{\mathchoice

{\setbox\ipbox=\hbox{$\displaystyle \left\langle\mathstrut
#1\right.$}

\vrule height\ht\ipbox width0.25pt depth\dp\ipbox}

{\setbox\ipbox=\hbox{$\textstyle \left\langle\mathstrut
#1\right.$}

\vrule height\ht\ipbox width0.25pt depth\dp\ipbox}

{\setbox\ipbox=\hbox{$\scriptstyle \left\langle\mathstrut
#1\right.$}

\vrule height\ht\ipbox width0.25pt depth\dp\ipbox}

{\setbox\ipbox=\hbox{$\scriptscriptstyle \left\langle\mathstrut
#1\right.$}

\vrule height\ht\ipbox width0.25pt depth\dp\ipbox}

}\right. }

\newcommand{\dirack}[1]{\left. \mathrel{\mathchoice

{\setbox\ipbox=\hbox{$\displaystyle \left.\mathstrut
#1\right\rangle$}

\vrule height\ht\ipbox width0.25pt depth\dp\ipbox}

{\setbox\ipbox=\hbox{$\textstyle \left.\mathstrut
#1\right\rangle$}

\vrule height\ht\ipbox width0.25pt depth\dp\ipbox}

{\setbox\ipbox=\hbox{$\scriptstyle \left.\mathstrut
#1\right\rangle$}

\vrule height\ht\ipbox width0.25pt depth\dp\ipbox}

{\setbox\ipbox=\hbox{$\scriptscriptstyle \left.\mathstrut
#1\right\rangle$}

\vrule height\ht\ipbox width0.25pt depth\dp\ipbox}

} #1\right\rangle}

\newcommand{\bz}{\mathbb{Z}}

\newcommand{\br}{\mathbb{R}}
\newcommand{\bc}{\mathbb{C}}

\newcommand{\bn}{\mathbb{N}}

\def\blfootnote{\xdef\@thefnmark{}\@footnotetext}

\renewcommand{\mod}{\operatorname{mod}}

\input xy
\xyoption{all}
\usepackage{amssymb}

\def\-{^{-1}}

\begin{document}
\title[Duality questions for operators, spectrum and measures]{Duality questions for operators, spectrum and measures}
\author{Dorin Ervin Dutkay}
\blfootnote{Work supported in part by the National science Foundation.}
\address{[Dorin Ervin Dutkay] University of Central Florida\\
	Department of Mathematics\\
	4000 Central Florida Blvd.\\
	P.O. Box 161364\\
	Orlando, FL 32816-1364\\
U.S.A.\\} \email{ddutkay@mail.ucf.edu}

\author{Palle E.T. Jorgensen}
\address{[Palle E.T. Jorgensen]University of Iowa\\
Department of Mathematics\\
14 MacLean Hall\\
Iowa City, IA 52242-1419\\}\email{jorgen@math.uiowa.edu}
\thanks{} 

\begin{abstract}
   We explore spectral duality in the context of measures in $\br^n$, starting with partial differential operators and Fuglede's question (1974) about the relationship between orthogonal bases of complex exponentials in $L^2(\Omega)$ and tiling properties of $\Omega$, then continuing with affine iterated function systems. We review results in the literature from 1974 up to the present, and we relate them to a general framework for spectral duality for pairs of Borel measures in $\br^n$, formulated first by Jorgensen and Pedersen.
\end{abstract}
\maketitle \tableofcontents
\section{Introduction}\label{intr}

We shall address a number of questions in the spirit of Mark Kac's question  ``Can you hear the shape of a drum?'' Since the question is open ended, it has already received a rich variety of answers in the literature, see for example \cite{BC86,Lap91,GWW92,Lap93,LP93,BCDS94,LP96,LF00,Sun85}. However the point of view we take here is different in at least three ways: while the instances alluded to in these references make a link between geometry and spectrum with the use of an appropriate Laplace operator, we will instead here rely on suitable generalizations of Fourier transform/Fourier series-ideas. While there is a rich and varied set of possibilities for geometries in the literature, we will here focus primarily on those that are based on two notions: tiling of some kind, and on matrix-scaling. While earlier transform tools have focused on one variable, for example the spectral variable of a Laplace operator, we will concentrate here on multi-variable geometries that can be accounted for  by means of self-similarity, for example that which is dictated by affine iterated function systems (IFS).

     A central motivation from the early days of functional analysis was the idea of using the Spectral Theorem in the analysis of partial differential operators (PDOs). By ``the Spectral Theorem'', we here refer to von Neumann's version, expressing a possibly unbounded selfadjoint operator as an integral of a spectral resolution, i.e., a projection-valued (Borel) measure on the real line. So, in going from one variable to several, say from one to $n$,  we will be looking for $n$ commuting families of projection-valued measures; or equivalently for a single projection valued measure defined on the Borel subsets in $\br^n$.

    But in applications to boundary value problems, one typically does not have self-adjointness. Rather one is faced with a system of formally selfadjoint operators (Hermitian operators with a common dense domain) in a suitable Hilbert space. The simplest PDE problem of this kind is obtained for a bounded open set $\Omega$ in $\br^n$, by taking a partial derivative operator in each coordinate direction. For a common domain for these $n$  operators we may then take differentiable functions on $\Omega$ which vanish on the boundary. These operators may individually have selfadjoint extensions. Hence for each of them, we will get a projection valued measure in one variable, but unless there is a system of strongly commuting selfadjoint extensions, there will not be a joint spectral resolution; i.e., there will not be a single projection valued measure on $\br^n$ that simultaneously diagonalizes the coordinate directions. However when commuting selfadjoint extensions exist, there will be $n$ associated unitary one-parameter groups; one for each coordinate direction. And if the selfadjoint extensions are strongly commuting, they then piece together to a unitary representation $U$ of the additive group $\br^n$ acting on $L^2(\Omega)$ where on $\Omega$ we use the restriction of Lebesgue measure from $\br^n$. One checks that locally the representation $U$ acts on $\Omega$ by translation. The choice of commuting selfadjoint extensions (if any) amounts of a functional analytic assignment of ``boundary conditions'' for $\Omega$. Hence the expectation is that under the global action by $U$, $\Omega$ will ``wrap around'' and retrace itself, and as a result one would expect that the joint spectrum, i.e., the spectrum of the representation $U$ will be periodic or quasi-periodic of some kind. This is born out with simple examples using Fourier series.

     Much of this program was made precise in a pioneering paper by Bent Fuglede \cite{Fug74} which in turn motivated a rich literature up to the present. Some of the early papers are \cite{Jor79,Jor80,Jor82a,Jor82b,Ped87,JoPe92}. In the time that followed these, the problem in fact took off in a variety of different new and fascinating directions, several of them motivated by the kind of selfsimilarity that is encoded by IFSs, including the question of building Fourier series and Fourier transforms on fractals.

We will review some of the results in the literature that stemmed out of Fuglede's work and were centered around Fuglede's conjecture. This is by no means a complete account of the subject, but merely a list of results and problems aimed mainly at students that might be interested in the subject. We would like to recommend also the web-site maintained by Izabella \L aba: {\sf http://www.math.ubc.ca/~ilaba/tiling.html}, which contains a very enjoyable exposition about spectral sets and tilings, open questions, and results.

We begin with a short review of Fuglede's paper \cite{Fug74}. In 1974, while working on commuting extensions of partial differential operators on domains in $\br^n$, Bent Fuglede proved the following theorem:
\begin{theorem}\label{thfug}
In $\br^n$ let $\Omega$ be a Nikodym region, i.e., a connected open set on which every distribution of finite Dirichlet integral is itself in $L^2(\Omega)$. The existence of $n$ commuting self-adjoint operators $H_1,\dots,H_n$ on $L^2(\Omega)$ such that each $H_j$ is the restriction of $-i\frac{\partial}{\partial x_j}$ (acting in the distribution sense) is equivalent to the existence of a set $\Lambda\subset\br^n$ such that the restrictions to $\Omega$ of the exponential functions $e_\lambda(x):=\exp(2\pi i\lambda\cdot x)$ form a total orthogonal family in $L^2(\Omega)$.
\end{theorem}

\begin{definition}\label{defsp}
When the set $\Lambda$ exists we say that $\Omega$ is a {\it spectral set} and $\Lambda$ is called a {\it spectrum }for $\Omega$.
\end{definition}

\begin{proof}[Some comments on the proof of Theorem \ref{thfug}] Suppose commuting self-adjoint extension operators $H_1,\dots, H_n$ exist. Let $E:\mathcal B(\br^n)\rightarrow\operatorname*{Proj}(L^2(\Omega))$ be the projection valued measure defined on the Borel sets in $\br^n$ and determined by
\begin{equation}\label{eqf1}
H_j=\int_{\br^n}\lambda_j\,dE(\lambda_1,\dots,\lambda_n),\quad j=1,\dots,n.
\end{equation}
Fuglede \cite{Fug74} proved that $E$ is necessarily atomic. Hence there is a discrete set $\Lambda\subset\br^n$ such that $\operatorname*{supp} E=\Lambda$. Moreover, Fuglede proved that for $\lambda\in\Lambda$, the projection $E(\{\lambda\})$ is one-dimensional and that 
\begin{equation}\label{eqf2}
E(\{\lambda\})=\bc e_\lambda
\end{equation}
where $e_\lambda(x):=e^{2\pi i\lambda\cdot x}$, $x\in\Omega$. 

From \eqref{eqf1}-\eqref{eqf2} one gets 
\begin{equation}\label{eqf3}
H_je_\lambda=\lambda_je_\lambda,\mbox{ for }\lambda=(\lambda_1,\dots,\lambda_n)\in\Lambda.
\end{equation}
Since $\sum_{\lambda\in\Lambda}E(\{\lambda\})=I_{L^2(\Omega)}$ one concludes that $\{e_\lambda\,|\,\lambda\in\Lambda\}$ is an orthogonal basis in $L^2(\Omega)$.

Note that, given \eqref{eqf1}-\eqref{eqf2}, the converse implication is simple. Having $\{e_\lambda\,|\,\lambda\in\Lambda\}$ an orthogonal basis allows one to define $\{H_j\}_j$ by \eqref{eqf3}.
\end{proof}

Fuglede proposed a way of checking whether such a spectrum exists in his celebrated conjecture:

\begin{conjecture}{\bf [Fuglede]}
A subset $\Omega$ of $\br^n$ of finite positive Lebesgue measure is spectral if and only if $\Omega$ tiles $\br^n$ by translations (up to measure zero), i.e., there exits a family of translations $\Gamma\subset\br^n$ such that $\Omega+\gamma$, $\gamma\in\Gamma$ is partition of $\br^n$ up to measure zero.
\end{conjecture}

\begin{definition}\label{deftil}
A measurable set $\Omega$ of finite positive Lebesgue measure is said to {\it tile} $\br^d$ by translations if there exists a set $\Gamma$ such that $\Omega+\gamma$, $\gamma\in\Gamma$ forms a partition of $\br^d$ up to Lebesgue measure zero. In this case $\Gamma$ is called a {\it tiling set} for $\Omega$.
\end{definition}

Fuglede proved the conjecture in the special case when the spectrum $\Lambda$ or the tiling set $\Gamma$ are full lattices in $\br^n$, but the question remained open until recently for the less structured cases. In 2004, Terence Tao gave a counterexample showing that one of the implications of Fuglede's conjecture is false in 5 and higher dimensions: if the dimension $n\geq5$ then there exists a finite domain in $\br^n$ that is spectral but does not tile $\br^n$ by translations. Subsequently, M{\'a}t{\'e} Matolcsi \cite{Mat05} found a counterexample for the same implication in dimension 4, and then with Mihail Kolountzakis for dimension 3 \cite{KM06}. For the reverse implication, tiles with no spectra were constructed for dimension $\geq5$ in \cite{KM06}, for dimension 4 in \cite{FR06}, and then for dimension 3 in \cite{FMM06}. 
Thus both implications of Fuglede's conjecture were disproved for dimensions higher than 3. As far as we know, at the time of writing this paper, the conjecture is still open in both directions in dimensions 1 and 2.

On the other hand, many results were obtained showing showing that the Fuglede conjecture is true for a more restrictive class of sets $\Omega$.
\begin{enumerate}
\item[$\bullet$]
A triangle is not spectral \cite{Fug74}.
\item[$\bullet$]
The unit ball is not spectral. Announced in \cite{Fug74}, proved in \cite{IKP99} and improved in \cite{Fug01} showing that the unit ball in $\br^n$ does not have an infinite family of pairwise orthogonal exponentials.
\item[$\bullet$]
If $\Omega$ is a bounded, open, convex set which is non-symmetric in the sense that $\Omega-x\neq x-\Omega$ for all $x\in\br^n$, then $\Omega$ is not spectral \cite{Kol00}.
\item[$\bullet$]
No symmetric convex body with smooth boundary can be spectral \cite{IKT01}. Moreover any set of orthogonal exponentials must be finite \cite{IR03}.
\item[$\bullet$]
The Fuglede conjecture holds in the special case where $\Omega$ is a convex compact set in the plane. More precisely, $\Omega$ admits a spectrum if and only if it is either a quadrilateral or a hexagon \cite{IKT03}.
\item[$\bullet$]
A union of two intervals of the form $(0,r)\cup(a,a+1-r)$ is spectral iff $r=1/2$ and $a=k/2$ for some $k\in\bz$ \cite{Lab01}.
\item[$\bullet$]
If $\Omega$ tiles $\br^+$ by translations then $\Omega$ tiles $\br$ by translations and is a spectral set \cite{PW01}.
\end{enumerate}

One of the complications in the Fuglede conjecture is the fact that a set can possess several spectra. To bypass this complication Lagarias and Wang introduced the concept of ``universal spectrum''. And since the Fuglede conjecture can be easily formulated in a more general context, for locally compact abelian groups, we include the definition here. More importantly, Fuglede's conjecture for $\br^n$ is intimately connected to the same conjecture for simpler groups, such as $\bz^n$ or $\bz_p$. The counterexamples constructed by Tao et. al. in high dimensions were based exactly on this observation, e.g. the construction of spectral sets in $\bz_{p_1}\times\dots\times\bz_{p_n}$ which do not tile this group can be used to obtain a non-tiling spectral set in $\br^n$, by adding the cube $[0,1)^n$.

\begin{definition}\label{defg}
 Let $G$ be a locally compact abelian group. Let $\Omega$ be a subset of $G$ of positive finite Haar measure. Then $\Omega$ is called a {\it spectral} set if there exists a family of characters $\Lambda$ in the Pontryagin dual $\hat G$ such that their restrictions to $\Omega$ form a total orthogonal family in $L^2(\Omega$). We say that $\Omega$ {\it tiles} $G$ if there exists a set $\Gamma\subset G$ such that $(\Omega+\gamma)_{\gamma\in\Gamma}$ forms a partition of $G$ up to Haar measure zero.
 \end{definition}
 
 \begin{definition}\label{defus}
 A {\it universal spectrum} \cite{LW97} for a set $\Gamma\subset\br^n$ is a set $\Lambda$ that is a spectrum simultaneously for {\it all} bounded measurable sets $\Omega$ that tile $\br^n$ with tiling set $\Gamma$.

 A {\it universal tiling} \cite{PW01} for a set $\Lambda\subset\br^n$ is a set $\Gamma$ that is a tiling simultaneously for {\it all} bounded measurable sets $\Omega$ that have $\Lambda$ as a spectrum.
 \end{definition}
 
 \begin{definition}\label{deffa}
 A {\it factorization} of a finite abelian group $G$, written $G=A\oplus B$, is one in which every $g\in G$ has a unique representation $g=a+b$ with $a\in A$ and $b\in B$.

 For $\mathcal A\subset \frac1{N_1}\bz\times\dots\times\frac{1}{N_n}\bz$ we say that a set $\mathcal B$ is a complementing set for $\mathcal A$ if $\mathcal B\subset\frac{1}{N_1}\bz\times\dots\times\frac{1}{N_n}\bz\cap[0,1]^n$ and $A:=\mathcal A(\mod\bz^n)$ and $B:=\mathcal B(\mod\bz^n)$ yield a factorization 
 $$A\oplus B=\bz_{N_1}\times\dots\times\bz_{N_n}:=\left(\frac{1}{N_1}\bz\times\dots\times\frac{1}{N_n}\bz\right)/\bz^n.$$

  Let $\bz_n^+=\{0,1,\dots,n-1\}$. We call $A\subset\bz^+$ a {\it direct summand} of $\bz_n^+$ if there exists a $B\subset\bz^+$ such that $A\oplus B=\bz_n^+$. We call a subset $\mathcal T$ a {\it strongly periodic set} if there exists an $n\in\bn$ and a direct summand $A\subset\bz^+$ of $\bz_n^+$ such that $\mathcal T=\alpha(A\oplus n\bz)$ for some non-zero $\alpha\in\bz$. See \cite{PW01}.
 \end{definition}
 
 \begin{conjecture}\label{cousc}\cite{LW97}
 {\bf [Universal Spectrum Conjecture]} Let $\Gamma:=\bz^n+\mathcal A$, where $\mathcal A\subset\frac{1}{N_1}\bz\times\dots\frac{1}{N_n}\bz$ such that $A:=\mathcal A(\mod\bz^n)$ admits some factorization $A\oplus B=\bz_{N_1}\times\dots\times\bz_{N_n}$. Then $\Gamma$ has a universal spectrum $\Lambda$ of the form $N_1\bz\times\dots\times N_2\bz+\mathcal L$, with $\mathcal L\subset\bz^n$.
 \end{conjecture}
 
 Several results are available in connection to the Universal Spectral Conjecture.
 \begin{enumerate}
 \item[$\bullet$]\cite{Fug74}
 Any lattice $\Gamma$ in $\br^d$ has a universal spectrum $\Lambda=\Gamma^*$ where $\Gamma^*$ is the dual lattice
 $$\Gamma^*=\{\lambda\in\br^d\,|\,\lambda\cdot \gamma\in\bz\mbox{ for all }\gamma\in\Gamma\}.$$
 \item[$\bullet$]\cite{LW97}
 Let $\Gamma=\bz^n+\mathcal A$ with $\mathcal A=\{a_1,\dots,a_m\}\subset\frac{1}{N_1}\bz\times\dots\times\frac{1}{N_n}\bz$. Then a set $\Lambda=(N_1\bz\times\dots \times N_n\bz)+\mathcal L$ with $\mathcal L\subset\bz$ is a universal spectrum for $\Gamma$ if and only if $\Lambda$ is a spectrum for each of the sets 
 $$\Omega_{\mathcal B}:=\left[0,\frac{1}{N_1}\right]\times\dots\times\left[0,\frac{1}{N_n}\right]+\mathcal B,$$
 for all complementing sets $\mathcal B$ for $\mathcal A$ in $\bz_{N_1}\times\dots\times\bz_{N_n}$. 
  \item[$\bullet$]\cite{LW97}
  If the group $\bz_N$ has the strong Tijdeman property, then any tile set $\Gamma=\bz+\frac{1}{N}\mathcal A$ with $\mathcal A\subset\bz$ has a universal spectrum $\Lambda=N\bz+\mathcal L$ for some $\mathcal L\subset\bz$.
   \item[$\bullet$]\cite{LW97}
   Let $\Gamma=\bz+\frac{1}{N}\mathcal A$, where $\mathcal A\subset\bz$ is such that $\mathcal A(\mod N)$ is a complementing set of $\bz_N$. If either $|\mathcal A|$ or $N/|\mathcal A|$ is a prime power, then $\Gamma$ has a universal spectrum $\Gamma=N\bz+\mathcal L$ with $\mathcal L\subset\bz$.
   \item[$\bullet$]\cite{PW01}
   Every strongly periodic set has a universal tiling and a universal spectrum which are strongly periodic.
   
 \end{enumerate}

 Particular attention has been given to sets of the form $B+[0,1)^n$ where $B\subset\br^n$ is finite. The next theorem explains the interest.
 
 \begin{theorem}\cite{KM06a, DJ08}
 Suppose $B\subset\bz^n$ is a finite set. Then $B$ is a spectral set in $\bz^n$ if and only if $B+[0,1)^n$ is a spectral set in $\br^d$. 
 \end{theorem}

 \begin{enumerate}
 \item[$\bullet$] A set $\Lambda$ is a spectrum for the cube $[0,1)^n$ if and only if $\Lambda$ is a tiling set for $[0,1)^n$ \cite{IP98,LRW00}. Conjectured and proved in dimension $\leq3$ in \cite{JP99}.
 \item[$\bullet$] Let $B\subset\bz$ be a finite set. If the polynomial $\sum_{b\in B}z^b$ is irreducible in $\bz[X]$ then the Fuglede conjecture holds for $B+[0,1)$ \cite{KL03}.
 \end{enumerate}
 
 Most examples of translational tiles exhibit some periodicity. The Periodic Tiling Conjecture (a topic which precedes Fuglede's paper) claims that periodicity must always be present, at least for regions of $\br^n$.
 
 \begin{definition}
 A set $\Gamma$ in $\br^n$ is called {\it periodic} if there exists a invertible matrix $R$ and a finite subset $L$ of $\br^n$ such that $\Gamma=L+R\bz^n$.
 \end{definition}
 
\begin{conjecture}\label{coptc}{\bf [Periodic Tiling Conjecture]}\cite{GS87}
Let $\Omega$ be a region in $\br^n$, i.e., a closed subset of $\br^n$ which is the closure of its interior, has finite positive Lebesgue measure, and has boundary $\partial \Omega$ of measure zero. If $\Omega$ tiles $\br^n$ by translations then it has a periodic tiling.
\end{conjecture}

\begin{enumerate}
\item[$\bullet$]
The one-dimensional case of this conjecture was proved by Lagarias and Wang in \cite{LW96}. 
\item[$\bullet$]
In dimension 2, the conjecture holds whenever $\Omega$ is a topological disk with smooth boundary \cite{GN89,Ken92}.
\item[$\bullet$]
The Periodic Tiling Conjecture holds for convex polytopes \cite{Ven54}.

\end{enumerate}

\section{Spectral measures}

While Fuglede's original version of the duality question was for the restriction of Lebesgue measure to bounded open subsets in $\br^n$, the early work suggested that spectral pairs automatically appear to have hidden selfsimilarity. The simplest way to make this precise is to follow Hutchinson (Theorem \ref{thhut} \cite{Hut81}) below; i.e., to consider measures that are fixed under a finite set of attractive affine mappings in $\br^n$. While Hutchinson considered more general iterated function systems (IFSs), we will restrict attention here to the affine class, see Definition \ref{defai}. This setting includes such important structures as Cantor sets and Sierpinski gaskets (a geometric and recursive method for creating gaskets: start with a triangle and cut out the middle piece. This results in three smaller triangles to which the process is continued. The nine resulting smaller triangles are cut in the same way, and so on, indefinitely); and even a more general framework with pairs of measures $(\mu, \nu)$ as introduced in section 4 below.

\begin{definition}\label{defsm}
Let $\mu$ be a Borel probability measure on $\br^n$. Then $\mu$ is called a {\it spectral measure} if there exists a subset $\Lambda$ of $\br^n$ such that the family of exponential functions $\exp(2\pi i\lambda\cdot x)$ with $\lambda\in\Lambda$ form an orthonormal basis for $L^2(\mu)$. In this case $\Lambda$ is called a {\it spectrum } of the measure $\mu$.
\end{definition}

 {\it Iterated function systems} (IFS) in $\br^d$ are natural generalizations of more familiar Cantor sets on the real line. Like their linear counterparts, they arise as limit sets $X$ for recursively defined dynamical systems. While the functions used may be affine, the limit $X$ itself will typically be a highly non-linear object, and will include complicated geometries. They arise in operator algebras and in representation theory; and they form models for ``attractors'' in dynamical systems arising in nature. For $d = 2$, the Sierpinski gasket is a notable example, and there is a variety of possibilities for $d>2$ as well. Each affine IFS $X$ possesses (normalized) invariant measures $\mu$, naturally associated with the system at hand (denote by $X$ the support of the measure). We focus here on the case when the maps are affine.
 
  Two approaches to IFSs have been popular: one based on a discrete version of the more familiar and classical second order Laplace differential operator of {\it potential theory}, see \cite{KSW01, Kig04, LNRG96};  and the second approach is based on {\it Fourier series}, which we discuss here. The first model is motivated by infinite discrete network of resistors, and the harmonic functions are defined by minimizing a global measure of resistance, but this approach does not rely on Fourier series. In contrast, the second approach begins with Fourier series, and it has its classical origins in lacunary Fourier series \cite{Kah86}.

 Here we make precise Hutchinson's equilibrium measures, i.e., measures $\mu$  that are fixed under a finite set of attractive affine mappings in $\br^n$. While Hutchinson considered more general iterated function systems (IFSs), we will restrict attention here to the affine class. 
\begin{definition}\label{defai}
Let $R$ be a $n\times n$ expansive matrix, i.e., all eigenvalues $\lambda$ have $|\lambda|>1$. Let $B$ be a finite subset of $\br^n$. Define the maps
$$\tau_b(x)=R^{-1}(x+b),\quad(x\in\br^n,b\in B).$$
We say that $(\tau_b)_{b\in B}$ is an {\it affine iterated function system} (affine IFS). We will denote by $N:=\#B$ the cardinality of $B$.
\end{definition}

\begin{theorem}\label{thhut}\cite{Hut81}
There exists a unique compact set $X_B$ such that 
\begin{equation}\label{eqxb}
X_B=\bigcup_{b\in B}\tau_b(X_B).
\end{equation}
The set $X_B$ is called the \textup{attractor} of the IFS $(\tau_b)_{b\in B}$.

Let $(p_b)_{b\in B}$ be a list of positive probabilities $\sum_{b\in B}p_b=1$.
There exists a unique Borel probability measure $\mu_{B,p}$ on $\br^n$ such that
\begin{equation}\label{eqmub}
\int f\,d\mu_{B,p}=\sum_{b\in B}p_b\int f\circ\tau_b\,d\mu_{B,p},
\end{equation}
for all continuous functions $f$ on $\br^n$. The measure $\mu_B$ is called \textup{the invariant measure} of the IFS. We denote by $\mu_B$, the measure $\mu_B$ in the case when all probabilities $p_b$ are equal $1/N$. The measure $\mu_{B,p}$ is supported on the attractor $X_B$.
\end{theorem}

 In understanding spectral duality for Hutchinson's equilibrium measures, we must establish orthogonality of a suitable set of Fourier frequencies. We will do this in a recursive manner. Just as the affine fractal and its equilibrium measure result by iterating the contractive affine mappings with scaling matrix $R^{-1}$, using expansion with $R^T$ (the transposed of $R$) we get an iteration in the large. We introduce certain cycles in the analysis of this fractal in the large, and the aim is to get the spectrum $\Lambda$ this way.

But to get started, we must initialize the Fourier duality. This is done in the definition below, we call this structure Hadamard triples. 
\begin{definition}
Let $R$ be a $n\times n$ an expansive integer matrix, i.e., all entries are integers. Let $B, L$ be subsets of $\bz^n$ of the same cardinality $\#B=\#L=:N$, and $0\in B$, $0\in L$. We say that $(R,B,L)$ is a Hadamard triple if the matrix
$$\frac{1}{\sqrt N}(e^{2\pi i R^{-1}b\cdot l})_{b\in B,l\in L}$$
is unitary.
\end{definition}

\begin{conjecture}\label{codj}\cite{DJ07}
If $(R,B,L)$ is a Hadamard triple then the measure $\mu_B$ is spectral.
\end{conjecture}

In finding the spectrum of a fractal measure, the following function $W_B$ plays a crucial role. Its heroes determine the existence and the formula for the spectrum.

\begin{definition}\label{defcy}
Let
\begin{equation}\label{eqwb}
W_B(x)=\left|\frac{1}{N}\sum_{b\in B} e^{2\pi ib\cdot x}\right|^2,\quad(x\in \br^n).
\end{equation}
For the finite set $L$ consider the IFS $\sigma_l(x)=(R^T)^{-1}(x+l)$, $x\in\br^n$, $l\in L$. A finite set $$C=\{x_0,\sigma_{l_1}x_0,\dots,\sigma_{l_{p-1}}\dots\sigma_{l_1}x_0\}$$ is called a {\it $W_B$-cycle} if there is some $l_p\in L$ such that 
$\sigma_{l_p}\dots\sigma_{l_1}x_0=x_0$ and $W_B(\sigma_{l_k}\dots\sigma_{l_1}x_0)=1$ for all $k\in\{1,\dots, p\}$.
\end{definition}

We now summarize some known results regarding iteration in the large: details about $W_B$-cycles (Definition \ref{defcy}) in the analysis of this ``fractal in the large'': this is our approach to Fourier spectrum $\Lambda$, and therefore to spectral pairs $(\mu, \Lambda)$.
\begin{enumerate}
\item[$\bullet$]\cite{JP98}
The conjecture holds for the Cantor measure $\mu_B$ obtained with $R=4$ and $B=\{0,2\}$. In this case a spectrum for the measure $\mu_B$ is $\Lambda=\{\sum_{k=0}^n4^k a_k\,|\, a_k\in\{0,1\}\}.$
\item[$\bullet$]\cite{Str98}
The conjecture holds if the zero set of the function $W_B$ is disjoint from the set $\sum_{k=1}^\infty (R^T)^{-k}L$.
\item[$\bullet$]\cite{LW02}
The conjecture holds in dimension $n=1$ if the only $W_B$-cycle is the trivial one $\{0\}$. 
\item[$\bullet$]\cite{DJ06}
The conjecture holds in dimension $n=1$, and a spectrum is the smallest set $\Lambda$ with the property that $\Lambda$ contains $-C$ for all $W_B$-cycles, and $R^T\Lambda+L\subset\Lambda$.
\item[$\bullet$]\cite{DJ07,DJ08a}
The conjecture holds in higher dimensions if a ``reducibility'' condition is satisfied.
\item[$\bullet$]\cite{CHR97}
The conjecture holds when $B$ is a complete set of representatives for $\bz^n/R\bz^n$. In this case the measure $\mu_B$ is a multiple of the Lebesgue measure on $X_B$, the set $X_B$ tiles $\br^n$ by a lattice hence the spectrum of $\mu_B$ is also a lattice. In this case $X_B$ is called a {\it self-affine tile}. See e.g.  \cite{GH94,LW96jlms,LW96s,LW97a,Wan99,LL07}.
\end{enumerate}

More general fractal spectral measures were constructed in \cite{Str00}. 
\begin{definition}
Let $\mathcal B$ and $\mathcal L$ be finite subsets of $\br^n$ of the same cardinality $N$. We say that $\{\mathcal B,\mathcal L\}$ is a {\it compatible pair} if the $N\times N$ matrix $\frac{1}{\sqrt N}(e^{2\pi ib\cdot l})_{b\in\mathcal B,l\in\mathcal L}$ is unitary. We denote by $\delta_{\mathcal B}$ the atomic measure
$$\delta_{\mathcal B}=\frac{1}{N}\sum_{b\in\mathcal B}\delta_b.$$
A {\it compatible tower} is a sequence of compatible pairs $\{\mathcal B_0,\mathcal L_0\}, \{\mathcal B_1,\mathcal L_1\},\dots$ with $B_j\subset M_j^{-1}\bz^n$ and $L_j\subset\bz^n$, and matrices $R_j\in\operatorname*{GL}(n,M_{j-1}\bz^n)$ for $j\geq 1$.
\end{definition}

\begin{theorem}\label{thstr}\cite{Str00}
Given an infinite compatible tower whose compatible pairs and expanding matrices are all chosen from a finite set of compatible pairs and expanding matrices, with $0\in L_k$ for all $k$, the infinite convolution product measure
$$\mu=\delta_{\mathcal B_0}\ast(\delta_{\mathcal B_1}\circ R_1)\ast\dots\ast(\delta_{\mathcal B_k}\circ (R_k\dots R_1))\ast\dots$$
exists as a weak limit and is a compactly supported probability measure. The set of functions $e_\lambda$ with $\lambda\in \Lambda$
$$\Lambda:=L_0+R_1^TL_1+\dots+R_1^T\dots R_k^TL_K+\dots$$
is an orthonormal set in $L^2(\mu)$.

If the zero set of the trigonometric polynomial $\hat\delta_{\mathcal B_k}$ is separated from the set 
$$(R_k^T)^{-1}\dots(R_1^T)^{-1}L_0+(R_k^T)^{-1}\dots(R_2^T)^{-1}L_1+\dots+(R_k^T)^{-1}L_{k-1}$$
by a distance $\delta>0$, uniformly in $k$ for all large $k$, then $\Lambda$ is a spectrum for $\mu$.

\end{theorem}

Strichartz also constructed sampling formulas for functions with fractal spectrum \cite{Str00,HS01}, and proved some surprising results about the convergence properties of Fourier series on fractals \cite{Str06}: they are much better than their well known classical counterparts on $[0,1]$.

Verifying that an invariant measure for an affine IFS is {\it not} a spectral measure can be a difficult task too. Some results are available in this direction: in \cite{JP99} it is proved that the triadic Cantor set does not have more than 2 pairwise orthogonal exponentials. In \cite{Li07pems} conditions are given for a given natural candidate set not to be a spectrum. In \cite{LW06} it is proved that a class of absolutely continuous measures with good decay of the Fourier transform are not spectral, and that if the diameter of the support set for absolutely continuous measure is not much larger than the measure then the set tiles (and is spectral) by a lattice, and the measure is a multiple of the Lebesgue measure on that set. In \cite{DJ07mz} a class of Sierpinski-like fractals are analyzed and shown to be spectral or non-spectral. In \cite{Li07} conditions are given involving the prime factorization of $\det R$ that guarantee that the measure $\mu_B$ does not have certain spectra.

All known examples of spectral measures generated by affine iterated function systems use equal probabilities for the maps $\tau_b$. \L aba and Wang proposed the following conjecture in dimension 1: 
\begin{problem}\label{prob2}\textup{\cite{LW02}}
Let $\mu_{\mathcal B,p}$ be the invariant probability measure associated to the IFS $\tau_b(x)=\rho(x+\beta)$, $\beta\in \mathcal B$, and probabilities $(p_i)_{i=1}^N$, with $|\rho|<1$ and $\mathcal B\subset \br$, $\#\mathcal B=N$, $p_i>0$. 

Suppose $\mu_{\mathcal B,p}$ is a spectral measure. Then

\begin{enumerate}
\item $\rho=\frac1P$ for some $P\in\bz$.
\item $p_1=\dots=p_N=\frac1N$.
\item Suppose $0\in \mathcal B$. Then $\mathcal B=\alpha B$ for some $\alpha\neq 0$ and some $B\subset\bz$ such that there exists $L$ that makes $(P,B,L)$ a Hadamard triple.
\end{enumerate}
\end{problem}

A counterexample to part (iii) was found in \cite{DJ08}; however, this example is not ``far'' from satisfying the conditions in (iii), so with a slight modification, part (iii) of the \L aba-Wang conjecture might still be true.

The next problem, even though it does not involve spectral theory, is crucial for the analysis of affine IFS measures: it asks if there can be essential overlap.
\begin{problem}\label{prob3}
Suppose $R$ is an expansive integer matrix and $B\subset \bz$ is such that no two elements in $B$ are congruent $\mod R$, i.e., $(B-B)\cap R\bz^n=\{0\}$. Then the measure $\mu_B$ has no overlap, i.e., 
$$\mu_B(\tau_b(X_B)\cap\tau_{b'}(X_B))=0,\quad(b,b'\in B,b\neq b')$$
\end{problem}

As explained above, it is known \cite{JP98} that the triadic Cantor set measure does not have more than 2 pairwise orthogonal exponentials. But how about frames?
\begin{problem}
Are there any tight frames of exponentials on the triadic Cantor set? Construct frames of exponentials on the triadic Cantor set. 
\end{problem}

Recall that a {\it frame} for a Hilbert space $H$ is a set of vectors $(e_i)_{i\in I}$ such that there exist constants $A,B>0$ such that 
$$A\|f\|^2\leq\sum_{i\in I}|\ip{f}{e_i}|^2\leq B\|f\|^2,\quad(f\in H).$$
If $A=B$, the frame is called {\it tight}.

In this problem we are referring to the Hilbert space $L^2(\mu)$, where $\mu$ is the canonical Hausdorff (or equivalently, affine IFS) measure on the Cantor set.

\section{Finite spectral sets}
In \cite{JP95}, Fuglede's conjecture for finite unions of intervals of length 1 was proved to be equivalent to the following (see also \cite{KL03}):
\begin{conjecture}\label{cotz}
Let $A$ be a finite subset of $\bz$. The following conditions are equivalent:
\begin{enumerate}
\item[(T)]
$A$ tiles $\bz$ by translations;
\item[(S)]
$A$ is a spectral set.
\end{enumerate}
\end{conjecture}

Tiling $\bz$ with a finite tile is an old subject. An important result was obtained by Coven and Meyerowitz in \cite{CM99}. We describe it here:

Let
$$A(x)=\sum_{a\in A}x^a.$$

Let $\Phi_s(x)$ denote the $s$-th cyclotomic polynomial, defined inductively by
$$x^n-1=\prod_{s|n}\Phi_s(x).$$
We define $S_A$ to be the set of prime powers $p^\alpha$ such that $\Phi_{p^\alpha}(x)$ divides $A(x)$.
\begin{conjecture}\label{cocm}\cite{CM99}
$A$ tiles $\bz$ by translations if and only if the following conditions are satisfied:
\begin{enumerate}
\item[(T1)] $A(1)=\prod_{s\in S_A}\Phi_s(1)$.
\item[(T2)] If $s_1,\dots,s_k\in S_A$ are powers of different primes, then $\Phi_{s_1\dots s_k}(x)$ divides $A(x)$.
\end{enumerate}
\end{conjecture}

It was proved in \cite{CM99} that (T) implies (T1)-(T2), that (T) implies (T1), and that (T) implies (T2) under the additional assumption that $\#A$ has at most two distinct prime factors. In \cite{Lab02} it is proved that (T1)-(T2) implies that $A$ is spectral. Other results in connection to these problems include:

\begin{enumerate}
\item[$\bullet$]\cite{KL03} Conjectures \ref{cotz} and \ref{cocm} are true when $A(x)$ is an irreducible polynomial.
\item[$\bullet$]\cite{KL03} For polynomials of the form
$$\frac{x^{mn}-1}{x^m-1}$$
if $\#A=2$ then Conjectures \ref{cotz} and \ref{cocm} are true, and if $\#A\geq3$ then Conjecture $\ref{cocm}$ is true.
\item[$\bullet$]\cite{New77} If $A$ tiles $\bz$ by translations then all tilings are periodic, i.e., if $A\oplus B=\bz$ then $B=C\oplus M\bz$ for some $M\in\bn$, $C\in\bz$, with $\#A\cdot\#C=M$. In other words $A\oplus C=\bz_M(\mod M)$.
\item[$\bullet$]\cite{PW01,DJ08} All direct summands of $\bz_n^+$ (as in Definition \ref{deffa}) are spectral.
\end{enumerate}

\section{A general duality}
All these examples of spectral sets and measures can be incorporated in a more general duality theory proposed by Jorgensen and Pedersen in \cite{JP99}.

  The general framework for spectral duality is the study of spectral theory for special Hilbert spaces $L^2(\mu)$ for finite Borel measures on $\br^n$; typically measures $\mu$ with a particular support, e.g., an open subset in $\br^n$;  a tile under $\br^n$-translation; or an affine IFS-fractal.  In our consideration in Theorem \ref{thfug}, the measure $\mu$ is the restriction of Lebesgue measure to a Borel subset $\Omega$ of finite positive Lebesgue measure in $\br^n$.  In Fuglede's  original context for spectral pairs, the reasoning in formulas \eqref{eqf1}-\eqref{eqf3} above shows that if $\Omega$ is of spectral type, i.e., there is a set $\Lambda$ giving a Fourier basis in $L^2(\Omega)$, then the set $\Lambda$ is necessarily discrete. Moreover the dual measure is then the counting measure supported on $\Lambda$.

Further, our consideration of affine IFSs in section 3 also implies discreteness of the sets $\Lambda$ that allow possible spectral transforms. Based on earlier examples, the authors of \cite{JP99} suggested instead a more symmetric formulation of the notion of spectral pair: they ask for two measures $\mu$ and $\nu$ such that there is an isometric Fourier transform from $L^2(\mu)$  to $L^2( \nu)$, see \eqref{eqn1}, and mapping onto $L^2(\nu)$. This will then be perfectly symmetric, so that $(\mu, \nu)$ is a spectral pair in this sense if and only if $(\nu,  \mu)$ is a spectral pair. Furthermore this framework includes as a special case the usual Fourier transform in $L^2(\br^n)$.  So the latter is an instance where both sides in the duality have continuous spectrum.

     The generalization was motivated by two facts: We look for a version of  spectral pairs $(\mu, \nu)$ that applies to fractals in such a way that one side in the duality is a fractal in the small and the second a fractal in the large. By ``fractal in the large'' we want to include limits of graphs, see Theorem \ref{thstr} above. Such limits typically are fractals in their own right. The second motivation for \cite{JP99} was to extend Heisenberg's uncertainty principle for traditional Fourier duality to the more general version involving spectral pairs, see \cite{JP99} for detailed results.
\begin{definition}\label{def1}
        {\bf Fourier duality.} \cite{JP99}. Let $d\in\bn$, and let $\mu,\nu$ be a pair of Borel measures on $\br^d$. Consider Borel functions $f$ on $\br^d$ such that 
        \begin{equation}\label{eqn1}
        \widehat{f\,d\mu}(y):=\int_{\br^d}e^{-ix\cdot y}f(x)\,d\mu(x)\mbox{ is well defined}.
        \end{equation}
        
        We say that $(\mu,\nu)$ is a {\it spectral pair} if 
        $$\{\widehat{f\,d\mu}\,|\, f\in L^2(\mu)\}=L^2(\nu),$$
        and
        $$\int_{\br^d}|\widehat{f\,d\mu}(y)|^2\,d\nu(y)=\int_{\br^d}|f(x)|^2\,d\mu(x),\quad(f\in L^2(\mu)).$$
        
        \end{definition}
        
        Of course $\br^n$ can be replaced by a locally compact abelian group $G$, and the pair of measures $(\mu,\nu)$ will be supported on the group $G$ and its Pontryagin dual $\hat G$ respectively.

        \begin{problem}\label{prb1} {\bf Spectral pairs of measures in dynamics.} {\it (a) What are the spectral pair measures $(\mu,\nu)$ associated with affine iterated function systems (IFSs)?

 (b) The same question for Julia sets of a rational iteration $z\mapsto r(z)$, or $z\mapsto z^2+c$ ?

  Sample questions: If $\mu$ is an equilibrium measure supported on a particular fractal, and if $\nu$ is a second measure such that  $(\mu, \nu)$ is a spectral pair, must $\nu$ then be pure point mass? If not, what are the possible spectral types for $\nu$?

How do a priori {\it uncertainty relations} of Heisenberg type from classical harmonic analysis carry over to detailed conclusions for spectral pairs of measures?}
 \end{problem}

\begin{acknowledgements}
The authors are grateful for discussions over the years
with colleagues: Bent Fuglede, Deguang Han, Keri Kornelson, Jeff Lagarias, Ka-Sing Lau,
Steen Pedersen, Gabriel Picioroaga, Karen Shuman, Qiyu Sun, and Yang Wang. Further we are grateful to the
organizers of a Banff workshop in 2007, especially to Sergei Silvestrov.
\end{acknowledgements}

\bibliographystyle{alpha}	
\bibliography{prbl}

\end{document}